\documentclass[a4paper,12pt]{amsart}

\usepackage{amssymb}
\usepackage{amscd}

\newtheorem{theorem}{Theorem}
\newtheorem{corollary}{Corollary}[section]
\newtheorem{lemma}{Lemma}[section]
\newtheorem{prop}{Proposition}[section]
\newtheorem{definition}{Definition}[section]

\newtheorem{question}{Question}[section]

\def\E{\operatorname{E}}
\def\F{\operatorname{F}}
\def\G{\operatorname{G}}
\def\U{\operatorname{U}}
\def\SU{\operatorname{SU}}
\def\O{\operatorname{O}}
\def\SO{\operatorname{SO}}
\def\Sp{\operatorname{Sp}}
\def\GL{\operatorname{GL}}
\def\SL{\operatorname{SL}}

\def\Ad{\operatorname{Ad}}

\def\Aut{\operatorname{Aut}}

\def\diag{\operatorname{diag}}

\def\dim{\operatorname{dim}}

\def\exp{\operatorname{exp}}
\def\Exp{\operatorname{Exp}}

\def\Hom{\operatorname{Hom}}

\def\Int{\operatorname{Int}}

\def\Lie{\operatorname{Lie}}

\def\ord{\operatorname{ord}}
\def\rank{\operatorname{rank}}
\def\span{\operatorname{span}}

\newcommand{\fra}{\mathfrak{a}}
\newcommand{\frb}{\mathfrak{b}}
\newcommand{\frc}{\mathfrak{c}}

\newcommand{\fre}{\mathfrak{e}}
\newcommand{\frf}{\mathfrak{f}}
\newcommand{\frg}{\mathfrak{g}}
\newcommand{\frh}{\mathfrak{h}}

\newcommand{\frk}{\mathfrak{k}}

\newcommand{\frp}{\mathfrak{p}}

\newcommand{\frr}{\mathfrak{r}}
\newcommand{\frs}{\mathfrak{s}}
\newcommand{\frt}{\mathfrak{t}}
\newcommand{\fru}{\mathfrak{u}}

\newcommand{\bbC}{\mathbb{C}}

\newcommand{\bbR}{\mathbb{R}}
\newcommand{\bbZ}{\mathbb{Z}}

\newcommand{\lra}{\longrightarrow}

\begin{document}
\title{A note on closed subgroups of compact Lie groups}

\author{Jun Yu}
\address{Department of Mathematics, ETH Z\"urich, 8092 Z\"urich, Switzerland}
\email{jun.yu@math.ethz.ch}

\abstract{We study two notions ``Lie primitivity'' and ``strongly control fusion'' (``SCF") 
for closed subgroups of Lie groups (both are defined by Griess (\cite{Griess})), 
generalize a theorem of Borel-Serre and classify simple symmetric pairs which are also 
``SCF'' pairs.} 
\endabstract

%\date{December 2009}


\maketitle

\section{Introduction}\label{S: Introduction} 

%Subgroups of compact (or complex) Lie groups are systematically studied since Dynkin 
%(\cite{Dynkin}) in 1950s. In the mean time, Borel and other Mathematicians stated the 
%study of topology of homogeneous spaces (\cite{Borel2}). It is Borel who first found 
%that, torsions of cohomology of a compact Lie groups have relations with finite subgroups 
%(precisely,  elementary abelian $p$ subgroups) of them (\cite{Serre2}).  Since them, many 
%mathematicians studied the topology of homogeneous spaces or general $G$-spaces and 
%continuous or finite subgroups of Lie groups. We can't give a complete list of literature 
%here, but recommend you to take a look at the papers \cite{Andersen-Grodal-Moller-Viruel}, 
%\cite{Serre2}, \cite{Griess-Ryba}, \cite{Liebeck-Seitz} and references therein.    

A finite subgroup $F$ of a connected Lie group $G$ is called {\it Lie primitive} if whenever $H$ is 
a closed Lie subgroup such that $F\subset H\subset G$, then $H$ is finite or $H=G$ (\cite{Griess}). 
In the case $G=\SL_{n}$ (or $\SU(n)$), this is a stronger condition than linear primitivity. 
We extend Griess' definition to any closed subgroup $H$ of a compact Lie group $G$, it is called (Lie) 
{\it primitive} if $(Z_{G_0})_{0}\subset H$ and there is no closed subgroup $K$ such that  
$H\subset K\subset G$, $K_{0}\not= G_{0}$, and $H_{0}$ is a proper normal subgroup of $K_{0}$. 
It is called (Lie) {\it quasi primitive} if there exists a sequence 
$H_{0}\supset H_{1}\supset...\supset H_{s}$ of closed subgroups of $G$ such that $H_{0}=G$, 
$H_{s}=H$ and each $H_{i+1}$ is primitive in $H_{i}$. With this notion of Lie primitivity, we have 
a generalization of theorem of Borel-Serre (\cite{Serre}).  

\begin{theorem} \label{T:group reduction}
For a closed subgroup $S$ of a compact Lie group $G$, there exists
a quasi-primitive subgroup $H$ of $G$ with $H_{0}\not=1$ and $S
\subset H$, such that either $(1)$ or $(2)$ of the following holds,
\begin{itemize}
\item[(1)]{$H_{0}$ is a Cartan subgroup of $G$.}
\item[(2)]{$S A_{H}/A_{H}$ is Lie primitive in $H/A_{H}$.}
\end{itemize}
\end{theorem} 

The proof of this theorem and how it leads to the original Borel-Serre theorem is presented in 
Section \ref{S:B-S}.   

\smallskip

A closed subgroup $H$ of a group $G$ is said to ``strongly controls fusion'' (write as ``SCF'' for brevity) 
in $G$ if for any two subsets $B_{1}, B_{2}$ of $H$ which are conjugate in $G$, whenever 
$g\in G$ satisfies $gB_{1}g^{-1}=B_{2}$, there exists $h\in H$ such that
$gbg^{-1}=hbh^{-1}$ for any $b \in B_{1}$. This condition is very strong, it reduces the conjugacy question 
of subgroups of $H$ completely to that of $G$. In \cite{Griess}, it was showed the natural pairs 
$(\GL(n,\bbC),\O(n,\bbC))$, $(\GL(2n,\bbC),\Sp(n,\bbC))$, $(\SO(7,\bbC),\G_2(\bbC))$, 
$(\E_6(\bbC),\F_{4}(\bbC))$ are ``SCF'' pairs.  From results in \cite{Griess}, \cite{Larsen}, 
we also know that there are no ``SCF'' pairs $(G,H)$ with $G$ a classical Lie group and 
$H$ a compact (or complex) simple Lie group of type $\F_4,\E_6,\E_7,\E_8$. 

For a compact symmetric space $G/K$, it is known that $G=KBK$, where $B=\exp(\frb_0)$ and $\frb_0$ is a 
maximal abelian subspace of $\frp_0=(\frk_0)^{+}$ (cf, \cite{Helgason}). For a compact Lie group $G$ with 
a bi-invariant Riemannian metric and two closed subgroups $H,K$, we show $G=H\exp(\frp_0)K$, where $\frp_0$ 
is the orthogonal complement of the Lie algebras of $H,K$ in $\Lie G$. This kind of double coset decomposition 
enables us to prove simply the compact symmetric pairs $(\SO(7),\G_2)$, $(\E_6,\F_{4})$ are ``SCF'' pairs 
(which is implied by Griess' results). Moreover we show that, besides those known ones, there are no more 
compact simple symmetric pairs which are also ``SCF'' pairs (Theorem \ref{T:SCF symmetric pair}). The question 
of looking for more ``SCF'' pairs $(G,H)$ with $G$ a compact connected simple Lie group remains interesting. 

\begin{theorem}\label{T:SCF symmetric pair}
For a simple compact symmetric pair $(\frg_0, \frh_0)$ with $\frg_0$
semi-simple, if there exists an ``SCF" pair $(G, H)$ such that
$(\Lie G, \Lie H)=(\frg_0, \frh_0)$. Then $(\frg_0,\frh_0)$ falls
into the following list: $(\frk_0\oplus \frk_0,\Delta(\frk_0))$ with
$\frk_0$ compact simple, $(\mathfrak{su}(n),\mathfrak{so}(n))$,
$(\mathfrak{su}(2n),\mathfrak{sp}(n))$,
$(\mathfrak{so}(2n),\mathfrak{so}(2n-1))$, $(\fre_6, \frf_4)$.

If we require $H$ to be connected, then $(\frg_0,\frh_0)$ falls into
the following list: $(\frk_0\oplus \frk_0,\Delta(\frk_0))$ with
$\frk_0$ compact simple,
$(\mathfrak{su}(2n+1),\mathfrak{so}(2n+1))$,
$(\mathfrak{su}(2n),\mathfrak{sp}(n))$, \\
$(\mathfrak{so}(2n),\mathfrak{so}(2n-1))$, $(\fre_6, \frf_4)$.
\end{theorem}

{\it Notations.}
We fix the following notations. $G$ is a real Lie group, not necessary connected or  
semisimple, topology of $G$ is the Euclidean topology. 
$G_{0}$ is the connected component of $G$ containing identity 
element, and $\frg_0=\Lie G$ is the Lie algebra of $G$. A subalgebra 
$\frh_0\subset \frg_0=\Lie G$ is called a closed subalgebra 
if the Lie subgroup $\exp(\frh_0)\subset G$ is a closed subgroup of $G$. 

For a subset $X\subset G$, let $C_{G}(X)=\{g\in G| gxg^{-1}=x,\forall x\in X\}$ be  
the centralizer of $X$ in $G$, and $N_{G}(X)=\{g\in G| gxg^{-1}\in X, \forall x\in X\}$ 
be the normalizer of $X$ in $G$. $Z_{G}$ the center of $G$, $[G,G]$ is the commutator subgroup, 
$[\frg_0,\frg_0]$ is the commutator subalgebra. 

The groups considered in this note are mostly compact Lie groups, but those considered in 
\cite{Griess} were linear algebraic groups over an algebraic closed field of characteristic 
0. These two categories of groups are equivalent (cf, Proposition \ref{P:non-comapact strong}).    

%The results in this note may have applications to the study of
%harmonic analysis on $L^{2}(G/H)$ for $G$ a compact simple Lie group
%and $H$ a finite subgroup or a symmetric subgroup.

\medskip

\noindent{\it Acknowledgement.} The author would like to thank Professor David Vogan for teaching 
him Steinberg 's theorem and thank Professor Xu-Hua He for several suggestions on the writing of 
this note. This note is mostly inspired from Serre 's Lecture notes 
\cite{Serre} and Griess' paper \cite{Griess}. The author 's research is supported by a grant 
from Swiss National Science Foundation (Schweizerischer Nationalfonds).

\section{Primitive subgroups and a theorem of Borel-Serre}\label{S:B-S}

\subsection{Lie primitive subgroups}\label{SS: Primitive} 

In \cite{Griess}, Griess defined a notion of ``Lie primitivity'' for finite subgroups 
of linear algebraic groups as an analogue of ``linear primitivity'' for linear representations. 
Recall that, linear representation $V$ of a group $S$ is {\it primitive} if there are no 
nontrivial decompositions $V=\oplus_{i}V_{i}$ such that $S$ permutes $V_{i}$.  A finite 
subgroup $F$ of a connected Lie group $G$ is called {\it Lie primitive} if whenever $H$ is a  
closed linear algebraic subgroup such that $F\subset H\subset G$, then $H$ is finite or $H=G$.  

From the above definition, if a group $G$ with non-abelian Lie algebra has a finite 
(Lie) primitive subgroup $F$, $G$ must be semisimple, since otherwise 
$F\subset H=F\cdot Z_{G}\subset G$ with $H$ is non-finite and $H\neq G$. When $G=\SL(V)$, 
if $F$ is a finite Lie primitive subgroup, then the representation $V$ of $S$ must be 
linear primitive. On the other hand, primitive finite subgroups of $\SL(V)$ are not 
necessary linear primitive. For example, $A_5$ has an inclusion $A_5\subset\SO(3,\bbC)$, 
the corresponding inclusion $A_5\subset\SL(3,\bbC)$ is not Lie primitive by definition. 
But, as the corresponding representation of $A_5$ on $\bbC^{3}$ is irreducible and $A_5$ is 
a simple group, we see it is linear primitive.

We generalize Lie primitivity to all closed subgroups. First, we call a closed subalgebra 
$\frh_0$ of a Lie algebra $\frg_0$ {\it Lie primitive} if $N_{\frg_0}(\frh_0)=\frh_0$. 
This is also called self-normalizing in literature.

\begin{prop}\label{P:basic normalizer}
For a compact semi-simple Lie algebra $\fru_0$ and a closed subalgebra
$\frh_0$, we have $N_{\fru_0}(\frh_0)=N_{\fru_0}(N_{\fru_0}(\frh_0))$; moreover,   
$N_{\fru_0}(\frh_0)=\frh_0$ if and only if $C_{\fru_0}(\frh_0)\subset\frh_0$.
\end{prop}

\begin{proof} Note that, for a compact Lie algebra $\fru_0$ with two closed subalgebras 
$\frr_0,\frs_0$, if $\frs_0$ is an ideal of $\fru_0$ and $\frr_0$ is an ideal of $\frs_0$, 
then $\frr_0$ is also an ideal of $\fru_0$ as $\fru_0$ is a direct sum of its simple ideals 
and its center. Thus $\frh_0$ is an ideal of $N_{\fru_0}(N_{\fru_0}(\frh_0))$, which implies
$N_{\fru_0}(N_{\fru_0}(\frh_0))\subset N_{\fru_0}(\frh_0)$. Then
$N_{\fru_0}(N_{\fru_0}(\frh_0))= N_{\fru_0}(\frh_0)$.

$N_{\fru_0}(\frh_0)$ is of the form $N_{\fru_0}(\frh_0)=\frh_0\oplus\frs_0$ for another ideal 
$\frs_0$ of $N_{\fru_0}(\frh_0)$. Thus $N_{\fru_0}(\frh_0)=\frh_0$ if and only if $\frs_0=0$, 
which is equivalent to $C_{\fru_0}(\frh_0)\subset\frh_0$.
\end{proof}

\begin{corollary}\label{C:algebra reduction}
For a closed Lie subalgebra $\frh_0$ of a compact semi-simple Lie algebra $\fru_0$, let 
$\fra_0=Z_{\frh_0}$ and $\frs_0=[\frh_0,\frh_0]$. Then the following conditions are equivalent,
\begin{itemize}
\item[(1)]{$\frh_0$ is Lie primitive in $\fru_0$.}
\item[(2)]{$C_{\fru_0}(\frh_0)\subset\frh_0$.}
\item[(3)]{$\fra_0=Z_{C_{\fru_0}(\fra_0)}$, and $\frs_0$ is primitive in
$[C_{\fru_0}(\fra_0),C_{\fru_0}(\fra_0)]$.}
\item[(4)]{$\fra_0$ is a Cartan subalgebra of $C_{\fru_0}(\frs_0)$.}
\end{itemize}
\end{corollary}
\begin{proof}
This follows from Proposition \ref{P:basic normalizer} immediately.  
\end{proof}

For a compact Lie group $G$, we call $C_{G}([G_{0},G_{0}])$ {\it the abelian part} 
of $G$ and let $A_{G}=C_{G}([G_{0},G_{0}])$ denote it. One can see that 
$(A_{G})_0=(Z_{G_0})_0$ is the connected component of the center of $G_0$, so it is 
abelian. But $A_{G}$ is itself not always abelian. The Lie algebra of $A_{G}$ is the 
center of $\frg_0=\Lie G$. Note that, $A_{G}$ is the kernel of 
the adjoint homomorphism $\pi: G\longrightarrow \Aut([\frg_0,\frg_0])$, 
where $\frg_0=\Lie G$ is the Lie algebra of $G$. Since $[\frg_0,\frg_0]$ is 
semisimple, we have \[\Int([\frg_0,\frg_0])\subset G/A_{G} \subset
\Aut([\frg_0,\frg_0]),\] from which we get $A_{G/A_{G}}=1$. 

\begin{definition}\label{D:primitive} 
For a compact Lie group $G$, a closed subgroup $H$ is called (Lie) {\it 
primitive} if $(A_{G})_{0}\subset H$ and there are no closed
subgroups $K$ with
\[H\subset K\subset G, K_{0}\not= G_{0} \textrm{ and } H_{0} \textrm{ 
is a proper normal subgroup of } K_{0}.\] $H$ is called (Lie) 
{\it quasi primitive} if there is a sequence 
$H_{0}\supset H_{1}\supset...\supset H_{s}$ of closed subgroups of $G$, 
such that $H_{0}=G, H_{s}=H$ and each $H_{i+1}$ is primitive in $H_{i}$.
\end{definition}

In the definition of ``primitivity", the condition $(A_{G})_0\subset H$ induces 
$H/(A_{G})_0 \subset G/(A_{G})_0$, and the latter $G/(A_{G})_0$ is a semi-simple group 
of adjoint type. When $G$ is a connected semisimple Lie group and $H$ is a finite 
subgroup, Lie primitivity defined above is equivalent to that defined in \cite{Griess}. 
Examples of Lie primitive groups include full rank subgroups,
maximal closed subgroups, and those $H \subset SU(n)$ with $H$
semi-simple and the corresponding representation of $H$ on
$\bbC^{n}$ is irreducible.

\begin{lemma}\label{L:group to algebra}
For a compact groups $H\subset G$, let $\frg_0=\Lie G$, $\frh_0=\Lie H$.

If $\frh_0$ is Lie primitive in $\frg_0$ (that is,
$N_{\frg_0}(\frh_0)=\frh_0$), then $H$ is Lie primitive in $G$.

If $H$ is Lie primitive in $G$, then
$N_{\frg_0}(\frh_0)=\frh_0 \textrm{ or } \frg_0$ (that is, 
$\frh_0$ is Lie primitive in $\frg_0$ or
$\frh_0$ is an ideal of $\frg_0$).
\end{lemma}

\begin{proof}
Suppose $H$ is not Lie primitive, let $K$ be a closed subgroup of
$G$ with $\dim K<\dim G$ and $H_0$ is a proper normal subgroup of
$K_0$. Let $\frk_0=\Lie K$. Then $\frh_0\subset\frk_0\subset
N_{\frg_0}(\frk_0)$, $\frk_0\not=\frh_0$ and $\frk_0\not=\frg_0$. So
$N_{\frg_0}(\frh_0)\not=\frh_0$. Thus $\frh_0$ is not Lie primitive.

Suppose $N_{\frg_0}(\frh_0)\not=\frh_0 $ and
$N_{\frg_0}(\frh_0)\not=\frh_0$, let $K=N_{G}(H_0)$. Then $H\subset
K\subset G$, and $\Lie K=N_{\frg_0}(\frh_0)$. Thus $\dim K<\dim G$
and $H_0$ is a proper normal subgroup of $K_0$, so $H$ is not Lie
primitive.
\end{proof}

\medskip

When $N_{\frg_0}(\frh_0)=\frg_0$, that is, $\frh_0$ is an ideal of 
$\frg_0$. Let $G'=N_{G}(H_0)/H_0$ and $H'=H/H_0$. Then 
$\Lie G'=\frg_0/\frh_0$ is semi-simple, $H'\subset G'$ and $H'$ is 
a finite group. It is also clear that $H$ is Lie primitive in $G$ 
if and only if $H'$ is Lie primitive in $G'$.

%When $G$ is semi-simple, the subgroup $A_{G}$ is finite, $H$ is Lie primitive 
%in $G$ if and only if $H'=HA_{G}/A_{G}$ is Lie primitive in $G'=G/A_{G}$.

\begin{lemma}\label{L:centralizer-normalizer} 
If $S$ is a finite Lie primitive subgroup of a compact semi-simple
Lie group $G$, then $C_{G}(S), N_{G}(S)$ are finite.
\end{lemma}

\begin{proof}
Suppose $C_{G}(S)$ or $N_{G}(S)$ is not finite, let $H_1=N_{G}(S)$. Then 
$S\subset H_1\subset G$. As $S$ is primitive, so $G_0\subset H_1$, which 
implies $[G_0,S]=1$.  Choose any connected closed proper subgroup $K$ of 
$G_0$, let $H=SK$. Then $S\subset H\subset G$, $H$ is not finite and 
$H_0\neq G_0$. This contradicts to $S$ is primitive.  
\end{proof}

\begin{prop}\label{P:finite Weyl}
Let $G$ be a compact Lie group and $H$ a Lie primitive subgroup in
$G$, then $N_{G}(H)/H$ is a finite group. If $H$ is semi-simple, 
then $G$ is semi-simple and $C_{G}(H)$ is finite. If $H$ is finite, then for 
any normal subgroup $A$ of $H$ not contained in $A_{G}$, $N_{G}(A)$ is finite.
\end{prop}

\begin{proof}
All statements follow from Lemmas \ref{L:group to algebra}, \ref{L:centralizer-normalizer}, 
or from the definition. 
\end{proof}

\subsection{A generalization of a theorem of Borel-Serre}

In 1950s, Borel and Serre proved the following theorem, which is very useful for the 
study of finite subgroups of Lie groups. 

\begin{theorem}[Borel-Serre, \cite{Serre}] \label{T:Borel-Serre}
For any compact Lie group $G$ with a supersovable finite subgroup $S$, there always exists a 
maximal torus $T$ of $G$ such that $S\subset N_{G}(T)$.
\end{theorem}

\begin{lemma} \label{L:primitive closure}
For a compact semi-simple Lie group $G$ and a closed subgroup $S$ with $\dim S< \dim G$ 
which is not primitive, there exists a closed primitive subgroup $H$ of $G$, such that 
$S\subset H$, $H_0 \neq G_0$, and $S_0$ is a proper normal subgroup of $H_0$. 
\end{lemma}

\begin{proof} Prove by induction on $\dim G$. Since $S$ is not primitive, there exists 
a closed subgroup $H'\subset G$, such that $S\subset H'$, $S_0$ is a proper normal
subgroup of $H'_0$, and $H'_0\neq G_0$. Let $\frg_0,\frh'_0,\frs_0$ be Lie algebras of 
$G,H',S$ respectively. 

When $\frh'_0$ is not normal in $\frg_0$, let $H=N_{G}(H'_0)$, then $S\subset H\subset G$. 
By Proposition \ref{P:basic normalizer} and Lemma \ref{L:group to algebra}, $H$ is Lie 
primitive, then $H$ satisfies the conclusion of the lemma. 

When $\frh'_0$ is normal in $\frg_0$, consider the inclusion $SH'_0/H'_0 \subset SG_0/H'_0$, 
then \[\dim(SG_0/H'_0)<\dim G.\] By induction, we can find a closed subgroup $H$ with 
$SH'_0\subset H\subset SG_0$, $\dim H/H'_0<\dim G/H'_0$, $H'_0$ is normal in $H$, 
$\dim H/H'_0>0$ and $H/H'_0 \subset SG_0/H'_0$ is a primitive inclusion. Then 
$S\subset H\subset G$, $\dim H<\dim G$, $H$ is Lie primitive in $G$ and $S_0$ is a proper 
normal subgroup of $H_0$. 
\end{proof}

\begin{lemma} \label{L:group reduction 1}
For a closed subgroup $S$ of a compact Lie group $G$, there exists
a quasi-primitive subgroup $H$ of $G$ with $S\subset H$ and $S_0\neq H_0$, such that either 
$(1)$ or $(2)$ of the following holds,
\begin{itemize}
\item[(1)]{$H_{0}$ is abelian.}
\item[(2)]{$S A_{H}/A_{H}$ is Lie primitive in $H/A_{H}.$}
\end{itemize}
\end{lemma}

\begin{proof}
Prove by induction on $\dim G$. As it is enough to consider $SA_{G}/A_{G}\subset G/A_{G}$, 
so we may assume that $G$ is semi-simple and of adjoint type.  

If $S \subset G$ is primitive. Let $H=G$, so $H$ satisfies condition $(2)$ in the lemma.

If $S\subset G$ is not primitive. By lemma \ref{L:primitive closure}, there exists a primitive 
closed subgroup $H'\subset G$ with $S\subset H'$, $H'_0\not=G_0$, and $S_0$ is a proper normal 
subgroup of $H'_0$. 

When $H'_0$ is abelian, let $H=H'$. Then $H$ satisfies $1)$ in the conclusion.

When $H'_0$ is not abelian, one has $(A_{H'})_0\not=H'_0$, so $H'/A_{H'}$ is a compact Lie group 
of positive dimension. Reduce the question to $SA_{H'}/A_{H'}\subset H'/A_{H'}.$

Since \[\dim H'/A_{H'} < \dim G,\] by induction we can find a closed subgroup $H$ of $G$ with 
$SA_{H'}\subset H \subset H'$, $H_0\neq(SA_{H'})_0,$ $H/A_{H'}\subset H'/A_{H'}$ is a quasi-primitive 
inclusion, and $(H/A_{H'})_0$ is abelian or $SA_{H}A_{H'}/A_{H}A_{H'} \subset H/A_{H}A_{H'}$ is 
a primitive inclusion. Then $S\subset H\subset G$ with $H_0\neq S_0$. When $(H/A_{H'})_0$ is abelian, 
$H_0$ is also abelian. When $SA_{H}A_{H'}/A_{H}A_{H'} \subset H/A_{H}A_{H'}$ is a primitive
inclusion, $A_{H'} \subset H$ implies $A_{H'} \subset A_{H}$, so $S A_{H}/A_{H}$ is Lie primitive 
in $H/A_{H}$. 
\end{proof}

\begin{proof}[Proof of Theorem \ref{T:group reduction}] Prove by induction on $\dim G$. 
Reduce to $SA_{G}/A_{G} \subset G/A_{G}$, we may assume that $G$ is
semi-simple and of adjoint type.

From lemma \ref{L:group reduction 1}, there exists a closed
quasi-primitive subgroup $H'$ of $G$ with $S\subset H'$,
$H'_0\not=S_0$ and $S A_{H'}/A_{H'}$ is Lie primitive in $H'/A_{H'}$
or $H'_{0}$ is abelian.

When $S A_{H'}/A_{H'}$ is Lie primitive in $H'/A_{H'}$, let $H=H'$,
then $H$ satisfies $2)$ in the conclusion.

When $H'_{0}$ is abelian, if $H_0$ is a Cartan subgroup. Let $H=H'$,
so $H$ satisfies $1)$.

If $H_0$ is not a Cartan subgroup. Then $N_{G}(H'_0)/H'_0$ is a
compact Lie group of positive dimension. Reduce to $S H'_{0}\leq
N_{G}(H'_0)/H'_0$, since \[\dim N_{G}(H'_0)/H'_0< \dim G,\] by
induction there exists a closed Lie subgroup $H$ of $G$ with $S
H'_0\subset H$, $S_0H'_0\not=H_0$, $H/H'_0\subset N_{G}(H'_0)/H'_0$ is
a quasi-primitive inclusion, and $H_0/H'_0$ is a Cartan subgroup of
$N_{G}(H'_0)/H'_0$ or $S A_{H}A_{H'}/A_{H}A_{H'}$ is Lie primitive
in $H/A_{H}A_{H'}$. Then argue as in the proof for lemma \ref{L:group
reduction 1}, we deduce that $H_{0}$ is a Cartan subgroup of $G$ or
$S A_{H}/A_{H}$ is primitive in $H/A_{H}$.
\end{proof}

\begin{proof}[Theorem \ref{T:group reduction} implies Theorem \ref{T:Borel-Serre}] 
By Theorem \ref{T:group reduction}, we only need to show
case $(2)$ in theorem \ref{T:group reduction} never happen for a
super-solvable subgroup $S$.

If case $(2)$ happens, we will get a primitive embedding from a
quotient of $S$ to some $H/A_{H}$. Let $\frh_0$ be the Lie algebra
of $H$. Since $\Int([\frh_0,\frh_0])\subset H/A_{H}\subset\Aut([\frh_0,
\frh_0])$, we may assume that $S$ itself has a primitive embedding to some
$G=\Aut(\fru_0)$ with $\fru_0$ a compact semi-simple Lie algebra. $S$ is
super solvable, so we can choose a prime order element $\sigma \in
S$ such that $\langle\sigma\rangle$ is normal in $S$. A classical
theorem of A. Borel said $\fru_0^{\sigma}\not= 0$. In another hand,
$\fru_0^{\sigma}\neq\fru_0$ since the automorphism $\sigma$ is
non-trivial. Let $K=N_{G}(\langle\sigma\rangle)$. Then $S\subset K
\subset G$, $K_0\not=1$ and $K_0\not=G_0$, so $S$ is not primitive 
in $G$. \end{proof}

\subsection{Finiteness}

\begin{prop} \label{P:primitive finite}
For a compact Lie group $G$, there are only finitely many 
conjugacy classes of Lie primitive closed subgroups of $G$.
\end{prop}

\begin{lemma} \label{L:group finite}
For a finite group $S$ and a compact semi-simple Lie Group $G$, the number 
of conjugacy classes of homomorphisms from $S$ to $G$ is finite.

For a complex semi-simple Lie algebra $\frh$ and another complex
semi-simple Lie algebra $\frg$, there are only finite many orbits
for the $\Int(\frg)$ conjugation action on $\Hom(\frh, \frg)$.
\end{lemma}

Lemma \ref{L:group finite} is a theorem of A. Weil (cf, \cite{Weil}). His method of proof  
is to establish a connection between infinitesimal deformations of homomorphism 
and Lie group (or Lie algebra) cohomology (actually $H^{1}$ only) and show 
some cohomology vanishes (e.g.,Whitehead 's first lemma). Some generalization of 
Weil 's theorem was considered by Jinpeng An and his coauthors (cf, \cite{An-Wang}).

\begin{proof}[Proof of Proposition \ref{P:primitive finite}] 
First of all, by Definition \ref{D:primitive}, $G$ has a Lie primitive finite subgroup 
implies $G$ is semisimple. By Lemmas \ref{L:group finite} and \ref{L:group to algebra}, 
we only need to show the possible isomorphism types of Lie primitive finite subgroups of 
any compact Lie group $G$ is finite. 

Embed $G$ faithfully into $GL(n, \bbC)$ for some positive integer
$n$ and let $S\subset G$ be a primitive finite subgroup. If $S$ has
a solvable minimal normal subgroup $A$. $A$ must be elementary
abelian, so $A\cong (C_{p})^{m}$ for a prime $p$ and a positive
integer $m$. $S \subset N_{G}(A)$ and $S$ is Lie primitive imply
that $N_{G}(A)$ is finite. By Borel-Serre theorem \ref{T:Borel-Serre}, 
$A$ is contained in some $N_{G}(T)$ for a 
maximal torus $T$ of $G$. If $p>n$, $A$ must be contained $T$, wchih
contradicts with $N_{G}(A)$ is finite, so $p<n$. Then $A\cap T$ is
contained in the $p$ torsion subgroup of $T$. Thus $m\leq \dim
T+\ord_{p}(|W|)$, where $\ord_{p}(k)$ is the order of maximal $p$
power dividing $k$ and $W$ is the Weyl group of $G$. Thus the
possibility of $A$ up to conjugation is finite, so the number of
conjugation classes of $S$ is also finite.

If $S$ has no solvable minimal normal subgroups. Choose a minimal
normal subgroup $B$ of $S$. It must be that $B\cong E^{k}$ for a
finite simple group $E$ and a positive integer $k$. Choose an
element $x$ of prime order $p$ in $E$, then $(C_{p})^{k} \subset
GL(n, \bbC)$. Argue as in the previous paragraph, from
$(C_{p})^{k}\subset GL(n,\bbC)$, we have $k\leq n+ord_{p}(n!)\leq
2n-1$. Since $E$ has a faithful representation of dimension $n$, for
each prime $p>n$, Sylow $p-$ subgroup of $E$ is abelian and
generated by at most $n$ elements, so $E$ has only finite many
isomorphism types. Then the number of isomorphism types of $B$ is
finite. $S\subset N_{G}(B)$ and $S$ is Lie primitive imply that
$N_{G}(B)$ is finite. Then by Lemma \ref{L:group finite}, the number
of conjugation classes of $S$ is finite.
\end{proof}

In \cite{Griess-Ryba}, Griess and Ryba classified quasisimple finite subgroups 
of exceptional compact simple Lie groups. The results are summarized in 
Table QE there, in which they also describe which subgroups are Lie primitive. 

A class of Lie primitive finite subgroups of a compact simple Lie group $G$ 
very different with quasisimple finite subgroups are those subgroups $S$ with a solvable 
normal subgroup $A$ not contained in $Z_{G}$. Since a solvable minimal normal subgroup $A$ 
is an elementary abelian $p$ subgroup. From the inclusion of $A$ in $G$, we know some 
information for the possible inclusions of $S$ in $G$.  

We would like to come back to the study of Lie primitive finite subgroups and general 
finite subgroups of compact simple Lie groups in future.

\section{``Strongly controlling fusion" pairs}

The notion of ``strongly control fusion'' was introduced by Griess in \cite{Griess}. 
For a pair of groups $H\subset G$, $H$ is said to ``strongly controls fusion in $G$'' 
if for any two subsets $B_{1}, B_{2}$ of $H$ which are conjugate in $G$, whenever 
$g\in G$ satisfies $gB_{1}g^{-1}=B_{2}$, there exists $h\in H$ such that
$gbg^{-1}=hbh^{-1}$ for any $b \in B_{1}$. The last condition is also equivalent to $g=hc$ 
for some $h \in H$ and $c \in C_{G}(B_{1})$. For simplicity, we will write ``SCF" for 
``strongly controls fusion". 

Note that, for $H\subset G'\subset G$, if $H$ is ``SCF" in 
$G$, then it is also ``SCF" in $G'$; if $H$ is ``SCF" in $G'$ and $G'$ is ``SCF" in $G$, 
then $H$ is ``SCF" in $G$. For a group $G$ with a closed subgroup $H$ and a normal closed
subgroup $N$, if $(G, H)$ is an ``SCF" pair, then $(G/N, H/(N\cap H))$ is also an ``SCF" 
pair. On the converse direction, if $N\cap H=1$ and $(G/N, H/(N\cap H))$ is an ``SCF",  
then $(G, H)$ is an ``SCF" pair. 

%Since for any $S_1, S_2 \subset G$ and $g\in G$ with
%$gS_1 g^{-1}=S_2$, $(G/N, H/(N\cap H))$ is an ``SCF" implies that $\exists h\in H$ such 
%that $gsg^{-1}(hsh^{-1})^{-1} \in N$ for any $s \in S_1$. On the other hand, 
%$gsg^{-1}, hsh^{-1} \in H$ and $H\cap N=1$. So $gsg^{-1}=hsh^{-1}, \forall s \in S_1$.

The notion of ``SCF" can be discussed on any class of groups, in this paper we consider 
real Lie groups.  An immediate consequence of the condition ``$H$ strongly controls its 
fusions in $G$" is the map \[\Hom(S, H)/H \longrightarrow \Hom(S, G)/G\] is injective.

\begin{lemma} \label{L:SCF criterion}
For groups $H\subset G$, $H$ strongly controls fusion in $G$ if
and only for any two closed subgroups $S_{1}, S_{2}$ of $H$ which are conjugate in $G$, whenever 
$g\in G$ satisfies $(S_{1})^{g}=S_{2}$, there exists $h\in H$
such that \[gsg^{-1}=hsh^{-1}, \forall s \in S_{1}.\] It is also equivalent to for any $g\in G$, 
$g \in HC_{G}(H\cap g^{-1}Hg).$
\end{lemma}

\begin{proof}
For the first statement, the condition in the original definition
implies that in this proposition is clear. If the condition in
this proposition holds, for a triple $(B_1,B_2,g)$ in the original
definition, let $S_{i}=\overline{\langle B_{i}\rangle}, i=1,2$, then
$gS_1 g^{-1}=S_2,$ thus there exists $h \in H$ such that
$s^{g}=s^{h}$ for any $s \in S_{1}$, in particular $b^{g}=b^{h}$ for
any $b \in B_{1}$. So the condition in the original definition holds.

For the second statement, we just need to use the fact for all $g\in
G$, $H\cap g^{-1}H g$ is the maximal subset $B_1$ of $H$ such that
$gB_1 g^{-1} \subset H$.
\end{proof}

If $g \in HC_{G}(H\cap g^{-1}Hg)$, then for any $g'=xgy^{-1}$ 
($x,y \in H$), $g'\in HC_{G}(H\cap g'^{-1}Hg')$. Thus the condition 
$g\in HC_{G}(H_{g})$ depends only on the double coset $HgH$ rather 
than on $g$ itself. By this, when the pair $(G,H)$ has a simple double 
cosets decomposition, it will be easier to check whether $(G, H)$ is 
an "SCF" pair.

\subsection{Examples of ``SCF'' pairs} \label{SCFT examples}

\begin{prop}\label{P:non-comapact strong}
For a real semi-simple connected Lie group $G$ with a Cartan
involution $\Theta$ and a maximal compact subgroup $K=G^{\Theta}$,
$K$ strongly controls fusion in $G$. 
\end{prop}

\begin{proof} Let $\frg_0=\frk_0+\frp_0$ be the corresponding Cartan
decomposition on Lie algebra Level. Choose a maximal abelian
subspace $\fra_0 \subset \frp_0$ and let $A=\exp(\fra_0)$. We have (\cite{Knapp})  
$G=KAK$, and the map $\exp: \frp_0\longrightarrow \Exp(\frp_0)$ is a 
diffeomorphism. If $(B_{1})^{g}=B_{2}$ with $B_{1}, B_{2} \subset  K$ and $g\in G$.
We may assume that $g\in A$, then
\begin{eqnarray*} &&\forall x \in B_{1}, g^{-1}xg=\Theta(gxg^{-1})=gxg^{-1}
\\&&\Longrightarrow \forall x\in B_1, g^{2}x=xg^{2}, \\&&
\Longrightarrow \forall x\in B_1, gx=xg, \end{eqnarray*} that is,
$g\in C_{G}(B_{1})$. So $K$ strongly controls fusion in $G$. 
\end{proof}

\begin{prop}
For any group $H$, $\Delta(H)=\{(x,x):x\in H\}$ strongly controls fusion in $G=H\times H$.
\end{prop}
\begin{proof}
This follows from: for any $g=(h_1,h_2)\in G$ and $h=(x,x)\in \Delta(H)$,
\[gxg^{-1}\in\Delta(H)\Leftrightarrow[h_1^{-1}h_2,x]=1\Leftrightarrow gxg^{-1}=g'xg'^{-1},\] 
where $g'=(h_2,h_2)$.  
\end{proof}

\begin{prop}[Griess, \cite{Griess} Theorem 2.3] \label{P:O(n) and Sp(n)}
For any $n\geq 1$ and the pair $(U,K)=(\U(n),\O(n))$ or $(U,K)=(\U(2n),\Sp(n))$, $K$ strongly
controls its fusions in $U$.
\end{prop}

\begin{prop} \label{P:strong criterion}
For $H=\Sp(n),\SU(n),\U(n),\SO(2n+1),\O(n)$, and any $G\supset H$, $H$ is ``SCF'' in 
$G$ if and only if ``$\forall x,y\in H$, $x\sim_{G} y$ implies $x\sim_{H} y$''.
\end{prop}

\begin{proof} For the if part, for two closed subgroups $S_1,S_2 \subset H$, if some 
$g\in G$ such that $g S_1 g^{-1}=S_2$, then $gsg^{-1}\sim_{G} s$ for any $s\in S_1$. By the 
assumption of the lemma, $gsg^{-1}\sim_{H} s$. By character theory and Proposition 
\ref{P:O(n) and Sp(n)}, there exists $h\in H$ such that $gsg^{-1}=hsh^{-1}$ for any 
$s\in S_1$, so $H$ strongly controls fusion in $G$.

For the only if part, let $S_1=\overline{\langle x\rangle},
S_2=\overline{\langle y\rangle}$. If $gxg^{-1}=y$ for some $g\in G$,
then $gS_1 g^{-1}=S_{2}$. Then there exists $h\in H$ such that
$gsg^{-1}=hsh^{-1}$ for any $s \in S_1$, so $hxh^{-1}=gxg^{-1}=y$.
\end{proof}

\begin{prop} \label{P:On SCF}
$\SO(2n-1)$ strongly controls its fusions in $\SO(2n)$.
\end{prop}

\begin{proof}
This follows from Proposition \ref{P:strong criterion}.   
\end{proof}

%\begin{proof}
%Let $A_{\theta}=\left(\begin{array}{cc} \cos\theta& \sin\theta \\ -\sin\theta&
%\cos\theta\end{array} \right)$. Then \[T=\{\left(\begin{array}{ccccc} A_{\theta_1}
%&&&& \\&\ddots&&&\\&&A_{\theta_{n-1}}&&\\&&&1&\\&&&&1\\\end{array}\right):
%0\leq\theta_1,...,\theta_{n-1} \leq 2\pi\}\] is a maximal torus of
%$SO(2n-1)$, and \[T'=\{\left(\begin{array}{cccc}A_{\theta_1}&&&\\&\ddots&&\\&&A_{\theta_{n-1}}&\\
%&&&A_{\theta_{n}}\\\end{array}\right):0\leq\theta_1,...,\theta_{n-1},\theta_{n}
%\leq 2\pi \}\] is a maximal torus of $SO(2n)$. The Weyl groups
%\begin{eqnarray*}&&W=W(SO(2n-1),T) \cong (C_{2})^{n-1} \rtimes
%S_{n-1},\\&& W'=W(SO(2n),T')\cong ((C_{2})^{n})'\rtimes S_{n},
%\end{eqnarray*} where $((C_{2})^{n})'=\{(\epsilon_1,\epsilon_2,...,
%\epsilon_{n})\in (C_{2})^{n}:\epsilon_1\epsilon_2\cdots\epsilon_{n}=1\}$.

%We know that for a compact Lie group $G$ with a maximal torus $S$,
%any element in $G$ is conjugate to an element in the maximal torus
%$S$, and for any two elements $x, y\in S$, $x, y$ are conjugate in
%$G$ if and only they are conjugate by some Weyl group element.

%For any $x, y\in T$, it is easy to see $x\sim_{W'} y $ if and only
%if $x\sim _{W} y,$ thus $SO(2n-1)$ strongly controls its fusion type
%in $SO(2n)$ by Lemma \ref{strong criterion}.
%\end{proof}

\begin{prop}[Griess, \cite{Griess} Theorem 1]\label{P:G2 SCF}
$\G_2$ strongly controls fusion in $\SO(7)$.
\end{prop}

\begin{prop}[Griess, \cite{Griess} Theorem 3]\label{P:F4 SCF}
$\F_4$  strongly controls fusion in $\E_6$.
\end{prop}

%Actually, in \cite{Griess} Griess proved that the pairs \[(\GL(n,k),\O(n,k)) \textrm{ and }   %
%(\GL(2n,k),\Sp(n,k))\] are ``SCF'' pairs for any algebraic closed field $k$. Apply to the 
%case $k=\bbC$, Griess' theorem implies Proposition \ref{P:O(n) and Sp(n)}, which is more 
%or less the same as the statement ``two real (or quaternion) representations are isomorphic 
%if and only if their complexfications are isomorphic''. 
%From Proposition \ref{P:O(n) and Sp(n)}, we also know that the pairs
%$(SU^{\pm{1}}(n),O(n))$ and $(SU(2n+1,SO(2n+1)))$ are ``SCF" pairs. Here 
%$SU^{\pm{1}}(n)$ is the group of $n\times n$ unitary matrices of determinant 
%$\pm{1}$.

Actually, in the paper \cite{Griess} it is showed 
\[(\GL(n,k),\O(n,k)), (\GL(2n,k),\Sp(n,k))\] are ``SCF '' pairs for any algebraic closed 
field $k$ of characteristic 0, and  \[(\SO(7,\bbC),\G_2(\bbC)), (\E_6(\bbC),\F_4(\bbC))\]
are ``SCF'' pairs. These statements are stronger than Propositions \ref{P:O(n) and Sp(n)}, 
\ref{P:G2 SCF} and \ref{P:F4 SCF} here. Propositions \ref{P:O(n) and Sp(n)}, \ref{P:G2 SCF} 
and \ref{P:F4 SCF} can be proved in the same way. That is, use the criterion in Lemma 
\ref{L:SCF criterion}, amounts to show \[g\in HC_{G}(H\cap g^{-1}Hg)\] for any $g\in G$. 
The difficult part of the proof is to calculate $H\cap g^{-1}Hg$. 

%\smallskip

We need to use several propositions in the proof.  
\begin{prop}[\cite{Helgason}]\label{P:Double decom1}
For a compact symmetic pair $G/K$, let $\frb_0$ be a maximal abelian subspace of 
$\frp_0=\frk_0^{+}$ and $B=\exp(\frb_0)$, then $G=KBK$.
\end{prop}

\begin{prop} \label{P:Double decom}
Let $G$ be a compact connected Lie group with two closed connected
subgroups $H, K\subset G$. Endow $G$ with a bi-invariant Riemannian
metric. Let $\frg_0$, $\frh_0$, $\frk_0$ be the Lie algebras of $G$,
$H$, $K$ respectively and $\frp_0=(\frh_0+\frk_0)^{\perp}$ be the
orthogonal complement of $\frh_0+\frk_0$ in $\frg_0$. Then
\[G=H\Exp(\frp_0)K.\]
\end{prop}

\begin{prop}[\cite{Knapp}] \label{P:connected centralizer}
For $G$ any compact connected Lie group and $S\subset G$ a torus(a
torus means a commutative connected subgroup), the centralizer
$C_{G}(S)$ is connected.
\end{prop}

\begin{prop}[Steinberg, \cite{Carter} Pages 93-95]\label{P:connected centralizer2} 
For $G$ any compact connected and simply connected Lie group and
$x\in G$ an element, the centralizer $G^{x}$ is connected.
\end{prop}

\begin{proof}[Proof of Proposition \ref{P:Double decom}]
For any $g\in G$, $H$ and $gK$ are two compact smooth sub-manifolds of $G$. Thus 
there exists a smooth shortest curve $\gamma: [0,1]\longrightarrow G$ connecting $H$ and 
$gK$, with $\gamma(0)\in H, \gamma(1)\in gK$. $\gamma$ must be a geodesic curve, let 
$\gamma(0)=h\in H,\gamma(1)=gk$ for some $h\in H, k\in K$. The equation of $\gamma$ can be 
written as $\gamma(t)=h\exp(tX)$ for some $X\in \frg_0$, then
\[\gamma(t)=h\exp(X)\exp((t-1)X)=\gamma(1)\exp((t-1)X)=gk\exp((t-1)X).\]

Since $\gamma$ is a shortest curve connecting connecting $H$ and
$gK$, one has that \[(L_{h})_{\ast}(X)\perp T_{h}(H) \textrm{ and }
(L_{gk})_{\ast}(X)\perp T_{gk}(gK).\]
Then $X\perp \frh_0$ and $X\perp \frk_0$. Thus $X\in \frp_0$.
Finally, \[gk=h \exp(X)\Longrightarrow g=h\exp(X) k^{-1}\in
H\exp(\frp_0)K.\] \end{proof}

\begin{proof}[Proof of Proposition \ref{P:Double decom1}]
When $H=K$ is a symmmetric subgroup, we have (\cite{Knapp}) 
$\frp_0=\Ad(K)(\frb_0)=\{\Ad(k)X:k\in K,X\in\frb_0\}$, so Proposition 
\ref{P:Double decom1} follows from Proposition \ref{P:Double decom}. 
\end{proof}

\smallskip 

\begin{proof}[Proof of Proposition \ref{P:O(n) and Sp(n)}] 
These two pairs are compact symmetric pairs, there
exists an involution $\sigma$ of $U$ such that $K=U^{\sigma}$ in
each case. Let $\fru=\frk_0+\frp_0$ be decomposition of the Lie
algebra $\fru=\Lie U$ into $+1, -1$ eigenspaces of $\sigma$. 
Note that, $\sigma=\tau$ for the pair $(\U(n),\O(n))$, and
$\sigma=\tau \Ad(\left(\begin{array}{cc} 0&I_{n}\\-I_{n}&0\\
\end{array}\right))$ for the pair $(\U(2n),\Sp(n))$.
Choose a maximal abelian subspace $\frb_0\subset \frp_0,$ let $B=\exp(\frb_0)$, 
then $U=KBK$ by Lemma \ref{P:Double decom1}. In the case of 
$(U,K)=(\U(n),\O(n))$, we may choose $B=\textrm{the set of diagonal\ matrices}$; 
In the case of $(U,K)=(\U(2n),\sp(n))$, we may choose 
\[B=\{\diag\{b',b'\}|b'=\diag\{x_1,x_2,\dots,x_{n}\},x_1x_2\cdots x_{n}=1\}.\]  
By Lemma \ref{L:SCF criterion}, we only need to show, for any $b\in B$, 
there exist $k\in K$ and $c\in C_{U}(K\cap b^{-1}Kb)$ such that $b=kc$. 
Let $\tau$ be the complex conjugation. 

For $x\in K$, $x\in K\cap b^{-1}Kb\Longleftrightarrow
b^{-1}xb=\sigma(bxb^{-1})=bxb^{-1}\Longleftrightarrow [x,b^{2}]=1$.
So $K\cap b^{-1}Kb=K^{b^{2}}$.

For $(U,K)=(\U(n),\O(n))$, $b\in B=\textrm{the set of diagonal\ matrices}$, let
\[b=\diag\{e^{\pi i t_1},e^{\pi i t_2},..., e^{\pi i t_n}\},
-1<t_{j}\leq 1. \]  Define 
$c=\diag\{e^{\pi i s_1},e^{\pi i s_2},..., e^{\pi i s_n}\}$ and 
$k=\diag\{r_1,r_2,...,r_n\}$, where $s_{j}=t_{j}$ if $-1< t_{j}\leq 0$,
$s_{j}=t_{j}-1$ if $0< t_{j}\leq 1$; $r_{j}=1$ if $-1< t_{j}\leq 0$,
$r_{j}=-1$ if $0<t_{j}\leq 1$ . Then $b=kc$ with $k\in K$ and  $c\in C_{U}
(K^{b^{2}})$. %In fact, $k ,c\in B$, $b^{2}=c^{2}$ and $k^{2}=1$.

For $(U,K)=(\U(2n),\Sp(n))$, the proof is similar. For $b=\diag\{b',b'\}\in B$, 
where $b'\in U(n)$ is a diagonal matrix. Let
\[b'=\diag\{e^{\pi i t_1},e^{\pi i t_2},...,e^{\pi i t_n}\}, 
-1<t_{j}\leq 1. \] Define $c'=\diag\{e^{\pi i s_1},e^{\pi i s_2},
..., e^{\pi i s_n}\}$ and $k'=\diag\{r_1,r_2,...,r_n\}$, where $s_{j}=t_{j}$ 
if $-1< t_{j}\leq 0$, $s_{j}=t_{j}-1$ if $0< t_{j}\leq 1$; $r_{j}=1$ if 
$-1< t_{j}\leq 0$, $r_{j}=-1$ if $0<t_{j}\leq 1$. Let $c=\diag\{c',c'\}$ and 
$k=\diag\{k',k'\}$. Then $b=kc$, $k\in K$ and $c\in C_{U}(K^{b^{1}})$.
\end{proof}

Recall that $\G_2$ has a unique seven dimensional irreducible representation,
which is a real representation. So we have an inclusion $\G_2\subset
\SO(7)$. In Lie algebra level, we have an orthogonal
decomposition (with respect to Killing form on $\mathfrak{so}(7)$)
$\mathfrak{so}(7)=\mathfrak{g}_2\oplus V$, where $V=\bbR^{7}$ is isomorphic 
to the seven dimensional real representation of $\G_2$.
Let $G=\SO(7)$, $H=\G_2$.

\begin{lemma} \label{L:G_2 stabilizer}
For any $0\neq v\in V$ and any $1\neq g\in \exp(V)$, the stabilizer(centralizer) 
for the adjoint action of $G_2$ \[G_2^{v}\cong G_2^{g}\cong SU(3).\]

For any $g\in \exp(V)$, \[H\cap g^{-1}Hg=H^{g^{3}}\]
\end{lemma}

\begin{proof}
For the first statement, one has $\dim(G_2v)\leq \dim V-1=6 $, so
\[\dim G_2^{v}=\dim G_2-\dim(G_2v)\geq 14-6=8.\] Then $G_2^{v}$ is
of rank two. Since $v\neq 0$, so $G_2^{v}\neq G_2$.  Any proper
connected subgroup of $G_2$ with dimension at least $8$ is conjugate
to $SU(3)\subset G_2$, so $(G_2^{v})_0\cong SU(3)$. Since $HgH \neq
G$, so $\dim HgH\leq \dim G-1=20$. Then \[\dim H^{g}=\dim H+\dim
gHg^{-1}-\dim HgH\geq 14+14-20=8.\]   Since $1\neq g$, so $H^{g}\neq
H$.   Then $(H^{g})_0$ is conjugate to $SU(3)\subset G_2\subset
SO(7)$.  Let $g=\exp(v),v\in V$. Then $H^{v}\subset H^{g}$ and
$H^{g}_0=H^{v}_0$.

The inclusion $SU(3)\subset SO(7)$ is in the well-known way
\[A+Bi\longmapsto \left(\begin{array}{ccc}A&B&  \\-B&A& \\&&1
\\\end{array}\right).
\] Assume $H^{g}_0=SU(3)$, then \[g\in C_{SO(7)}(SU(3))=
\{\left(\begin{array}{ccc}aI_3&bI_3&  \\-bI_3&aI_3& \\&&1
\\\end{array}\right): a^{2}+b^{2}=1\}.\]
Let $Z$ denote $C_{SO(7)}(SU(3))$, which is connected and of
dimension one. Since $\exp(V)\subset C_{G}(H^{v})$, so $Z=\exp(V)$.
One sees that \[Z\cap H=\{\left(\begin{array}{ccc}aI_3&bI_3&
\\-bI_3&aI_3& \\&&1\\\end{array}\right): (a+bi)^{3}=1\}\subset SU(3),\] and
for any $1\neq g\in Z=\exp(V)$, $SO(7)^{g}=U(3)=SU(3)Z$ with
$SU(3)\subset G_2=H$. Then \[H^{g}=H\cap SO(7)^{g}=H\cap
SU(3)Z=SU(3)(H\cap Z)=SU(3).\] So $H^{v}=H^{g}\cong SU(3)$.

For the second statement, let $L=SO(8)/\langle -I\rangle$ and
consider the inclusion $G_2\subset SO(7)\subset
L=SO(8)/\langle-I\rangle$, there exists an order three automorphism
$\theta$ of $L$ such that $L^{\theta}=G_2$. Then for any $x\in H$,
$x\in H\cap g^{-1}Hg$ if and only if $\theta(g^{-1}xg)=g^{-1}xg$,
which is also equivalent to $(\theta(g)g^{-1})x(\theta(g)g^{-1})^{-1}=x$. 
Then \[H\cap g^{-1}Hg=H^{\theta(g)g^{-1}}.\]

Under the adjoint action of $G_2$, let $U$ be the orthogonal complement of 
$\frg_2$ in $\mathfrak{so}(8)$. Then $\mathfrak{so}(8)=\mathfrak{g}_2\oplus U$, 
$\mathfrak{so}(7)=\mathfrak{g}_2\oplus V$ and $V\subset U$. Moreover, for any 
$v\in V$, $[\theta(v),v]=0$ and $\theta^{2}v+\theta v+v=0$. Then 
$[\theta g,g]=1$ and $(\theta^{2}g)(\theta g)g=1$.

For any $x\in H$, if $x\in H^{\theta(g)g^{-1}}$, then $(\theta g)g^{-1}$ and 
$g^{-1}(\theta g)^{-2}=(\theta^{2}g)(\theta g)^{-1}=\theta((\theta g)g^{-1})$
commute with $x=\theta x$. Thus $g^{3}=((\theta g)g^{-1})^{-2}
(g^{-1}(\theta g)^{-2})^{-1}$ commutes with $x$ and so $H^{\theta(g)g^{-1}}
\subset H^{g^{3}}$. It is clear that $H^{g}\subset H^{\theta(g)g^{-1}}$. 
So $G^{g}\subset H^{\theta(g)g^{-1}}\subset H^{g^{3}}$.

When $g^{3}\neq 1$, by the first statement, we have $H^{g^{3}}\cong
H^{g}\cong SU(3)$. So \[H^{\theta(g)g^{-1}}=H^{g^{3}}=H^{g}.\]

When $g^{3}=1$, by the proof for the first statement, we have
$g\in Z\cap H\in H$. So $H^{\theta(g)g^{-1}}=H^{g^{3}}=H$
\end{proof}

\begin{proof}[Proof of Proposition \ref{P:G2 SCF}]
Let $G=SO(7)$, $H=G_2$. By Proposition \ref{P:Double decom}, we have
$SO(7)=G_2\exp(V)G_2$. By the remark following Proposition \ref{L:SCF
criterion}, we only need to show \[\forall g\in\exp(V), g\in
HC_{G}(H\cap g^{-1}Hg).\]

From the proof for Lemma \ref{L:G_2 stabilizer}, when $g^{3}\neq 1$, we
have $H\cap g^{-1}Hg=H^{g^{3}}=H^{g}$, so \[g\in C_{G}(H^{g})\subset
HC_{G}(H\cap g^{-1}Hg);\] when $g^{3}=1$, we have 
$g\in H\subset HC_{G}(H\cap g^{-1}Hg)$.
\end{proof}

Recall that, there exists an outer involution $\theta$ of $\E_6$ such that 
$\E_6^{\theta}=\F_4$. In Lie algebra level, we have an orthogonal
decomposition $\mathfrak{e}_6=\mathfrak{f}_4\oplus V$, where $V=\bbR^{26}$ is
an real irreducible representation of $F_4$. Let 
\[G=\E_6,\ H=G^{\theta}=\F_4,\ \frg_0=\fre_6,\ \frh_0=\frg_0^{\theta}
=\frf_4.\] Let $\fra_0\subset V$ be a maximal abelian subspace, then $\dim\fra_0=2$, 
$C_{\frg_0}(\fra_0)=C_{\frh_0}(\fra_0)\oplus\fra_0$ and $C_{\frh_0}(\fra_0)\cong
\mathfrak{so}(8)$. Choose a Cartan subalgebra $\frs_0$ of $C_{\frh_0}(\fra_0)$, then 
$\frt_0=\frs_0+\fra_0$ is a Cartan subalgebra of $\frg_0$. One can choose co-root vectors 
$\{H'_{j}\subset i \frt_0:1\leq j\leq 6\}$ so that 
\begin{eqnarray*}&& \theta(H'_1)=H'_6, \theta(H'_2)=H'_2, \theta(H'_3)=H'_5,
 \\&& \theta(H'_4)=H'_4, \theta(H'_5)=H'_3, \theta(H'_6)=H'_1.
\end{eqnarray*} Then $\fra_0=\span\{i(H'_1-H'_6),i(H'_3-H'_5)\}$.
For any $g\in \exp(V)$ and any $v\in V$, $G^{v},G^{g}$ are connected by Lemma 
\ref{P:connected centralizer}. So $G^{v},G^{g}$ are determined by the roots 
annihilated by $v,g$.
For roots $\beta_1,\beta_2,...,\beta_{s}$ in a root system, let
$\langle\beta_1,\beta_2...,\beta_{s}\rangle$ denote the sub-root-system generated 
by $\beta_1,\beta_2,...,\beta_{s}$. Let $Z_{G}=Z_{E_6}$ be the center of $E_6$, 
which is a cyclic group of order three.

\begin{lemma} \label{L:F_4 stabilizer}
For any $g\in\exp(\fra_0)-Z_{G}$, there exists $v\in\fra_0$ such
that $g^{2}=\exp(2v)$ and $G^{g^{2}}=G^{v}$.
\end{lemma}

\begin{proof}
For $v'=\pi i(a(H'_1-H'_6)+b(H'_3-H'_5))\in \fra_0, a,b\in\bbR$ and
$g=\exp(v')$, $g\not\in Z_{G}$ is equivalent to say ``at most one
number among $\{a+b,2a-b,a-2b\}$ is an integer".

If none of $a+b,2a-b,a-2b$ is an integer, then the roots annialated
by $g^{2}$ are in \[\span\{\alpha_2,\alpha_4,\alpha_3+\alpha_4+\alpha_5,
\alpha_1+\alpha_3+\alpha_4+\alpha_5+\alpha_6\} (\textrm{of type } D_4).\] Let 
$v=v'=\pi i(a(H'_1-H'_6)+b(H'_3-H'_5))$. Then $g^{2}=\exp(2v)$ and 
$G^{g^{2}}=G^{v}=C_{G}(\fra_0)$.

If $a+b\in\bbZ$, then the roots annialated by $g^{2}$ are in
\[\span\{\alpha_3+\alpha_4+\alpha_5,\alpha_2,\alpha_4,\alpha_1+\alpha_3,\alpha_5+\alpha_6\}
(\textrm{of type } D_5),\] which are also annialzated by  $\pi i
(a(H'_1-H'_6-H'_3+H'_5)), \pi i ((a+1)(H'_1-H'_6-H'_3+H'_5))$. Let
$v$ be either of \[\pi i (a(H'_1-H'_6-H'_3+H'_5)), \pi i
((a+1)(H'_1-H'_6-H'_3+H'_5)).\] Then $g^{2}=\exp(2v)$ and
$G^{g^{2}}=G^{v}=C_{G}(\fra_0)$.

If $2a-b\in\bbZ$, then the roots annialated by $g^{2}$ are in
\[\span\{\alpha_4,\alpha_2,\alpha_3+\alpha_4+\alpha_5,\alpha_1,\alpha_6\}(\textrm{of type } D_5),\]
which are also annialzated by  $\pi i (a(H'_1-H'_6+2H'_3-2H'_5)),
\pi i ((a+1)(H'_1-H'_6+2H'_3-2H'_5))$. Let $v$ be either of \[\pi i
(a(H'_1-H'_6+2H'_3-2H'_5)), \pi i ((a+1)(H'_1-H'_6+2H'_3-2H'_5)).\]
Then $g^{2}=\exp(2v)$ and $G^{g^{2}}=G^{v}$.

If $a-2b\in\bbZ$, then the roots annialated by $g^{2}$ are in
\[\span\{\alpha_1+\alpha_3+\alpha_4+\alpha_5+\alpha_6,\alpha_2,\alpha_4,\alpha_3,\alpha_5\}(\textrm{of type } D_5),\]
which are also annialzated by  $\pi i (a(2H'_1-2H'_6+H'_3-H'_5)),
\pi i ((a+1)(2H'_1-2H'_6+H'_3-H'_5))$. Let $v$ be either of \[\pi i
(a(2H'_1-2H'_6+H'_3-H'_5)), \pi i ((a+1)(2H'_1-2H'_6+H'_3-H'_5)).\]
Then $g^{2}=\exp(2v)$ and $G^{g^{2}}=G^{v}$.
\end{proof}

\begin{proof}[Proof of Proposition \ref{P:F4 SCF}]
Let $G=E_6$, $H=F_4$ and $\theta$ be an involutive automorphism of
$E_6$ such that $G^{\theta}=H$.  Let $V=\bbR^{26}=\frg_0^{-\theta}$
be the orthogonal complement of $\frh_0$ in $\frg_0$ and $\fra_0$ be
a maximal abelian subspace in $V$.   Then $G=H\exp(\fra_0)H$ by Proposition 
\ref{P:Double decom1}. By the remark following Proposition \ref{L:SCF criterion}, 
we only need to show \[\forall g\in\exp(\fra_0), g\in HC_{G}(H\cap g^{-1}Hg).\]

Since $H=G^{\theta}$ and $\theta|_{\fra_0}=-1$, for any  $x\in H$,
$x\in H\cap g^{-1}Hg$ if and only if
$g^{-1}xg=\theta(gxg^{-1})=gxg^{-1}$, which is also equivalent to
$g^{2}xg^{-2}=x$.  Then \[\forall g\in\exp(\fra_0), H\cap
g^{-1}Hg=H^{g^{2}}.\]

Moreover, for any $x\in\exp(\fra_0)$, $x\in H\cap\exp(\fra_0)$ if
and only if $x^{-1}=\theta(x)=x$, which is also equivalent to
$x^{2}=1$. Then $H\cap\exp(\fra_0)=\{x\in\exp(\fra_0):x^{2}=1\}$.

For any $g\in\exp(\fra_0)-Z_{G}$, by Lemma \ref{L:F_4 stabilizer},
there exists $v\in \fra_0$ such that $g^{2}=\exp(2v)$ and
$G^{g^{2}}=G^{v}$. Let $h=g\exp(-v)\in\exp(\fra_0)$ and
$c=\exp(v)\in\exp(\fra_0)$, so $h^{2}=1$. Then $g=hc$, $h\in H$ and
$c\in C_{G}(H^{g^{2}})$. So $g\in HC_{G}(H\cap g^{-1}Hg)$.

For $g\in \exp(\fra_0)\cap Z_{G}$, we have $g\in C_{G}(H\cap
g^{-1}Hg)\subset HC_{G}(H\cap g^{-1}Hg)$.
\end{proof}

\subsection{Compact symmetric pairs}
Compact symmetric pairs are classified by Elie Cartan in 1920s \cite{Borel}, 
\cite{Knapp} and \cite{Helgason} contain nice descriptions for this classification. 
They are important in differential geometry and representation theory. 

\begin{lemma}\label{L:SCF}
For any compact Lie group $G$ with a closed subgroup $H$, if $\rank
H=\rank G$ and $\dim H<\dim G$, then $(G, H)$ is not an ``SCF" pair.
\end{lemma}

\begin{proof}
Assume $(G,H)$ is an ``SCF" pair. Let $G_0, H_0$ be the connected
components of $G, H$ containing identity element respectively.

Since $(G, H)$ is an ``SCF" pair, $(G'=G/A_{G}, H'H/(A_{G} \cap H))$
is also an ``SCF" pair. $\rank H=\rank G$ implies that the connected
component $(A_{G})_0\subset H$, $\dim H< \dim G$ and $\rank H=\rank
G$ implies that $G_0$ is not abelian, that is, $G'=G/A_{G}$ is a
compact Lie group of positive dimension. Thus the pair $(G', H')$ is
still a pair of compact Lie groups(of positive dimensions) satisfies
\begin{eqnarray*}\rank H'=\rank G'\ \textrm{and}\ \dim H'< \dim G'.
\end{eqnarray*} Since we can reduce the question
to the pair $(G', H')$ and $A_{G/A_{G}}=1$, we may assume that
$A_{G}=1$ for the original $(G, H)$ pair.

Now we assume $A_{G}=1$, then $G$ is a compact semi-simple Lie group
of adjoint type.

At first, it must be that $H_0=G_0 \cap H$. If $H_0\not=G_0\cap H$,
then there exists $x\in G_0\cap H-H_0$. Choose a maximal torus
$T\subset H_0$, since $\rank H=\rank G$, $T$ is also a maximal torus
of $G_0$, thus there exists $y\in T$ such that $x\sim_{G_0} y$.
$x\not \in H_0$, $y\in H_0$, thus $x\not \sim _{H} y$. This
contradicts with the condition $(G, H)$ is an ``SCF" pair.

$T$ is a maximal torus of $H_0$, it is also a maximal torus of
$G_0$, let $\frt_0,\frh_0,\frg_0$ be Lie algebras of $T, H, G$
respectively. $G$ is of adjoint type implies $C_{G}(T)=T$. Let
$N_{G}(T)$, $N_{G_0}(T)$, $N_{H}(T)$, $N_{H_0}(T)$ be normalizer of
$T$ in $G, G_0, H, H_0$ respectively, then \begin{eqnarray*}N_{H}(T)
\cap N_{G_0}(T)=N_{H_0}(T).\end{eqnarray*} Let \begin{eqnarray*}&&
W_{G}=N_{G}(T)/T,\ W_{G_0}=N_{G_0}(T)/T, \\&& W_{H}=N_{H}(T)/T,\
W_{H_0}=N_{H_0}(T)/T. \end{eqnarray*} Then $W_{H_0} \subset W_{H},
W_{G_0} \subset W_{G}$ and $W_{H} \cap W_{G_0}=W_{H_0}$. $G_0$ is
semi-simple and $\rank H< \rank G$ imply that $W_{H_0}
\not=W_{G_0}$, so $W_{H}\not=W_{G}$.

Now we choose an element $x\in T$, which is regular with respect to
the $W_{G}$ action on $T$, that is, \[\forall 1\not=w \in W_{G},
w(x)\not=x.\] Then for any $n\in N_{G}(T)-N_{H}(T)$, $x\sim_{G}
nxn^{-1}$ but $x\not \sim_{H} nxn^{-1}$.
\end{proof}

\begin{proof}[Proof of Theorem \ref{T:SCF symmetric pair}]
Simple compact symmetric pairs $(\frg_0,\frh_0)$ outside the first
list include (cf, \cite{Knapp} or \cite{Helgason}): pairs $(\frg_0,\frh_0)$ 
with $\rank \frh_0=\rank\frg_0$, $(\mathfrak{so}(2p+2q+2),\mathfrak{so}(2p+1)
\oplus\mathfrak{so}(2q+1))$ with $1\leq p\leq q$ and $p+q \geq 3$, 
$(\fre_6,\frc_4=\mathfrak{sp}(4))$.

If there exists an ``SCF" pair $(G, H)$ such that $(\Lie G, \Lie H)$
is one of such pairs. We may first assume $A_{G}=1$, then
$\Int(\frg_0)\subset G\subset \Aut(\frg_0)$ and $\exp(\frh_0)\subset
H \subset N_{\Aut(\frg_0)}(\exp(\frh_0))$, here $\exp: \frg_0\lra
\Aut(\frg_0)$ is the exponential map for $\Aut(\frg_0)$.

When $\rank \frh_0=\rank \frg_0$, $(G, H)$ is not an ``SCF" pair
follows from the last lemma.

When
$(\frg_0,\frh_0)=(\mathfrak{so}(2p+2q+2),\mathfrak{so}(2p+1)\oplus
\mathfrak{so}(2q+1))$ with $p\geq 1, q\geq 2$, choose any
$0<\theta<\pi,$ then for \begin{eqnarray*}X=\left(
\begin{array}{cccc} I_{2p-1}&&& \\ &A_{\theta}&& \\ & & I_{2q-1}&
\\&&&A_{\theta} \end{array} \right) \in H, Y=\left(
\begin{array}{cccc} I_{2p+1}&&& \\ &I_{2q-3}&& \\ & & A_{\theta}&
\\&&&A_{\theta} \end{array} \right)\in H,\end{eqnarray*} it is always that $X \sim_{G}
Y$ but $X\not \sim _{H} Y$.

When $(\frg_0,\frh_0)=(\fre_6,\frc_4=\mathfrak{sp}(4))$, we may
assume that \[\Int(\fre_6)\subset G\subset \Aut(\fre_6),
\exp(\mathfrak{sp}(4))\subset H\subset
N_{\Aut(\fre_6)}(\mathfrak{sp}(4)).\] Then
$H_0=\exp(\mathfrak{sp}(4))\cong Sp(4)/\langle -I\rangle$. Since
$\Aut(\mathfrak{sp}(4))=\Int(\mathfrak{sp}(4))$ and there exists an
outer involutive automorphism $\theta$ of $\fre_6$ such that
$\exp(\mathfrak{sp}(4))=\Aut(\fre_6)^{\theta}$, so
\[N_{\Aut(e_6)}(\exp(\mathfrak{sp}(4)))=\exp(\mathfrak{sp}(4)) \times
\langle \theta \rangle\] and $H_0=G_0\cap H$.  $G_0$ has only two
conjugation classes of involutions but $H_0$ has three conjugation
classes of involutions, thus $(G, H)$ is not an ``SCF" pair.

Simple compact symmetric pairs in the first list but not in the
second list include only $(\mathfrak{su}(2n),\mathfrak{so}(2n))$.

When $(\frg_0, \frh_0)=(\mathfrak{su}(2n),\mathfrak{so}(2n))$ and
$H$ is connected, we may assume that \[(G, H)=(SU(2n)/\langle
-I\rangle, SO(2n)/\langle -I \rangle).\] Let $X=e^{\frac{\pi
i}{n}}\left(\begin{array}{cc} -1&\\&I_{2n-1} \end{array} \right)$.
Then $\Ad(X)(H)=H$ and $\Ad(X)|_{H}$ is an outer automorphism of
$H$. Thus there is none $Y\in H$, such that
\[\Ad(X)|_{H}=\Ad(Y)|_{H}.\]

We have showed
\begin{eqnarray*}&& (SU^{\pm{1}}(n),O(n)),
(SU(2n),Sp(n)), (SU(2n+1),SO(2n+1)),\\&& (SO(2n), SO(2n-1))),
(E_6,F_4)\end{eqnarray*} are all ``SCF" pairs . Thus each of the
pairs $(\frg_0,\frh_0)$ in the two lists can be realized as $(\Lie
G, \Lie H)$ for an ``SCF" pair $(G, H)$.
\end{proof}

\subsection{More on ``SCF" pairs.}

Up to conjugation, any $U(1)\cong H \subset U(n)=G$ is of the form
\[H(a)=\{h(z)=\diag\{z^{a_1},z^{a_2},...,z^{a_{n}}\}| z\in U(1)
\}\] for some sequence $a=(a_1,a_2,...,a_{n})$ of non-negative
integers with $a_1\leq a_2\leq...\leq a_{n}$ and
$\gcd(a_1,a_2,...,a_{n})=1$. Since $H$ is abelian, $H$ strongly
controls its fusions in $G$ if and only for any $z,w\in U(1)$,
$z\not= w$ implies $h(z)\not \sim _{G} h(w)$, which means any two
different elements in $H$ are not conjugate in $G$.

\begin{lemma}
$H(a)$ doesn't strongly control its fusions in $U(n)$ if and only if
there exist integers $m>k>1$, such that
$ka=(ka_1,ka_2,...,ka_{n})\equiv a(mod m)$ as multi-sets.

For $a=(1,1,2,...,n-1)$, $H(a)$ strongly control fusion in
$U(n)$

For $a=(1,2,3,...,n)$, $H(a)$ doesn't strongly control fusion
in $U(n)$
\end{lemma}

\begin{proof} For the first statement, if there exist integers
$m>k>1$, such that $ka=(ka_1,ka_2,...,ka_{n})\equiv a(mod\ m)$, let
$z=e^{2\pi i/m},w=e^{2k\pi i/m}$, then $z\not=w$ and $h(z)\sim_{G}
h(w)$.

Now we suppose there is none pair $(m, k)$ with $1< k< m$ such that
$ka\equiv a(mod\ m)$ as multi-sets, and there exists $z, w\in U(1)$
with $h(z)\sim_{G} h(w), z\neq w$.

We first show that $z$ and $w$ are of finite orders. Since
$h(z)\sim_{G} h(w)$, $(z^{a_1},z^{a_2},...,z^{a_{n}})$ differs with
$(w^{a_1},w^{a_2},...,w^{a_{n}})$ only by a permutation. This means
we can permute $a=(a_1,a_2,...,a_{n})$ to $b=(b_1,b_2,...,b_{n})$
such that $z^{a_{j}}=w^{b_{j}}, \forall j, 1\leq j \leq n$. If there
exists $j\neq k$ such that $a_{j}b_{k}-a_{k}b_{j}\not=0$. \[z^{a_{j}}=
w^{b_{j}} \textrm{ and }z^{a_{k}}=w^{b_{k}}\Longrightarrow z^{a_{j}b_{k}-a_{k}b_{j}}=1.\] 
So $z$ is of finite order. If for any $j\neq k$,
$a_{j}b_{k}-a_{k}b_{j}=0$. Then $b=a$ or $b=-a$. $b=a$ and
$\gcd(a_1,a_2,...,a_{n})=1$ imply $z=w$, which contradicts to $z\neq
w$. $b=-a$ imply $-a$ differs with $a$ by a permutation. Then for
any $m\geq 3, k=m-1,$ $ka\equiv a(mod\ m)$ as multi-sets, which
contradicts to the assumption "there is none pair $(m, k)$ with $1<
k< m$ such that $ka\equiv a(mod\ m)$ as multi-sets".

Thus $z$ is of finite order. Similarly, $w$ is of finite order.
Since $\gcd(a_1,a_2,...,a_{n})=1$, $(z^{a_1},z^{a_2},...,z^{a_{n}})$
generate the subgroup $\langle z\rangle$ of order $o(z)$ in $U(1)$
and $(w^{a_1},w^{a_2},...,w^{a_{n}})$ generate the subgroup $\langle
w\rangle$ of order $o(w)$ in $U(1)$.
$(z^{a_1},z^{a_2},...,z^{a_{n}})$ differs with
$(w^{a_1},w^{a_2},...,w^{a_{n}})$ only by a permutation. Thus these
two cyclic groups are equal. Let $m=o(z)=o(w), w=z^{k}, m>k>1$, then
$ka=(ka_1,ka_2,...,ka_{n})\equiv a(mod\ m)$.

The latter two statements follow from the first directly.
\end{proof}

For the group $SO(3)$, for each non-negative integer $n$, $SO(3)$
has an irreducible representation $\rho_{n}$ of dimension $2n+1$,
these representations are actually over $\bbR$, so we have
\[\rho_{n}: SO(3)\longrightarrow SO(2n+1) \subset SU(2n+1).\] Let
$H_{n}$ denote the image of $\rho(SO(3))\subset SO(2n+1)$.  Let
$A_{\theta}=\left(
\begin{array}{cc} \cos(\theta)& \sin(\theta) \\ -\sin(\theta)&
\cos(\theta) \\ \end{array} \right)$, then $T=\{\left(
\begin{array}{cc} A_{\theta}&  \\  & 1 \\ \end{array} \right)|
0\leq \theta \leq 2\pi \}$ is a maximal torus of $SO(3)$. Under
certain normalization, \begin{eqnarray*}\rho_{n}\left(
\begin{array}{cc} A_{\theta}&  \\  & 1 \\ \end{array} \right)=\left(
\begin{array}{cccc} A_{n\theta}& & &
\\&...& &  \\ & &A_{\theta} & \\ & & &1
\end{array} \right).\end{eqnarray*}

\begin{lemma}
For any $n\geq 2$, $\rho_{n}(SO(3))=H_{n}\subset SO(2n+1)$ doesn't
strongly control fusion in $SO(2n+1)$.

$(\rho_{1}\oplus \rho_2)(SO(3))=H\subset SO(8)$ strongly controls
fusion in $SO(8)$.
\end{lemma}

\begin{proof} For the first statement, since $n\geq 2$,
the Euler number $\phi(2n+1)>2$, then there exists $1<k<2n$(for
example, $k=n$) such that $(k,2n+1)=1$. Then for $\theta=2\pi
/(2n+1)$, $A_{\theta}$ and $A_{k\theta}$ are not conjugate in
$SO(3)$, but $\rho_{n}(A_{\theta})$ and $\rho_{n}(A_{k\theta})$ are
conjugate in $SO(2n+1)$, thus $\rho_{n}(SO(3))=H_{n}\subset
SO(2n+1)$ doesn't strongly control its fusions in $SO(2n+1)$.

By Proposition \ref{P:strong criterion}, to prove the second statement, we
just need to show for $\theta, \phi$, if images of $A_{\theta}$ and
$A_{\phi}$ in $H$ are conjugate in $G$, then $A_{\theta}$ and
$A_{\phi}$ are conjugate in $SO(3)$, the latter means
$\phi=\theta(mod\ 2\pi)$ or $\phi=-\theta(mod\ 2\pi)$. This is an
elementary number theoretic or elementary combinatorial problem. We leave 
it them to the reader.
\end{proof}

Till now, we have showed some natural pairs are ``SCF'' pairs and more are not, it 
remains an interesting question to find more ``SCF'' pairs $(G,H)$ with $\rank H$ relatively 
small compared to $\rank$, in particular we ask the following question. 

\begin{question} For a sequence of positive integers $a=(a_1, a_2,...a_{s}),
a_1\leq a_2\leq .... \leq a_{s}$, let $n=2\sum_{1\leq j\leq
s}{a_{j}}+s$, for which sequences $a$,
\begin{eqnarray*}\bigoplus_{1\leq j\leq s} \rho_{a_{j}}: SO(3) \lra
SO(n)\end{eqnarray*} is an ``SCF" inclusion? 
\end{question}

If $H\subset G$ is an ``SCF'' pair, then the study of conjugacy questions in $H$ is 
completely reduced to that in $G$. This is particular useful when $H$ is an exceptional 
simple Lie group and $G$ is a classical Lie group. Unfortunately, by results in 
\cite{Griess}, \cite{Larsen}, any compact simple Lie group of type 
$\F_4,\E_6,\E_7,\E_8$ doesn't have such an ``SCF'' inclusion. 
In the above, we illustrated different ways to detect whether a pair $H\subset G$ is an 
``SCF'' pair and classified ``SCF'' pairs which are also symmetric pairs.  
It is still interesting question to get more ``SCF'' pairs $H\subset G$ with $G$ a compact 
simple Lie group. Propositions \ref{L:SCF criterion} and \ref{P:Double decom} 
would be helpful to prove some pairs $H\subset G$ are ``SCF'' when $G$ is not 
too large compared to $H$.

\end{document}